\documentclass[10pt,oneside]{jm}
\usepackage{amsfonts,amssymb,amsmath,amscd}

\newcommand{\abs}{\vskip 0.5em\noindent}
\newcommand{\Abs}{\paragraph{}\hspace{-1em}}
\newcommand{\AbsT}[1]{\paragraph{\hspace{-1em} #1}}
\newcommand{\AbsTT}[1]{\paragraph*{#1}}

\newcommand{\Pf}{\paragraph*{Proof.}}

\newcommand{\Ex}{\AbsT{Example.}}

\newcommand{\Thm}{\AbsT{Theorem.}}

\newcommand{\Cor}{\AbsT{Corollary.}}

\newcommand{\bfi}{\noindent {\bf i)} }

\newcommand{\bfii}{\noindent {\bf ii)} }

\newcommand{\bfa}{\noindent {\bf a)} }
\newcommand{\bfb}{\noindent {\bf b)} }
\newcommand{\bfc}{\noindent {\bf c)} }

\newcommand{\QED}{\hfill $\sharp$}

\newcommand{\N}{\mathbb N}
\newcommand{\Z}{\mathbb Z}
\newcommand{\Q}{\mathbb Q}

\newcommand{\C}{\mathbb C}

\newcommand{\bX}{\mathbb X}

\newcommand{\cC}{\mathcal C}

\newcommand{\cI}{\mathcal I}
\newcommand{\cJ}{\mathcal J}
\newcommand{\cK}{\mathcal K}
\newcommand{\cL}{\mathcal L}

\newcommand{\cM}{\mathcal M}

\newcommand{\cO}{\mathcal O}

\newcommand{\al}{\alpha}

\newcommand{\Gm}{\Gamma}

\newcommand{\eps}{\epsilon}
\newcommand{\ia}{\iota}

\newcommand{\ph}{\varphi}

\newcommand{\zt}{\zeta}

\newcommand{\GL}{\text{GL}}

\newcommand{\id}{\text{id}}
\newcommand{\im}{\text{im}}

\newcommand{\lcm}{\text{lcm}}

\newcommand{\modcat}[2]{\textbf{mod}_{#1}\text{-}#2}
\newcommand{\projcat}[2]{\textbf{proj}_{#1}\text{-}#2}

\newcommand{\tr}{\text{tr}}

\newcommand{\ld}{,\ldots\hskip0em ,}
\newcommand{\lr}[1]{\langle #1\rangle}

\newcommand{\un}[1]{\underline{#1}}

\newcommand{\wti}[1]{\widetilde{#1}}
\newcommand{\wh}[1]{\widehat{#1}}
\newcommand{\mt}{\mapsto}

\newcommand{\ra}{\rightarrow}

\newcommand{\cn}{\colon}
\newcommand{\dcup}{\stackrel{.}{\cup}}

\newcommand{\spmid}{\,\mid\,}
\newcommand{\spnmid}{\,\nmid\,}

\newcommand{\sseq}{\subseteq}
\newcommand{\smin}{\setminus}
\newcommand{\emp}{\emptyset}

\newcommand{\tm}{\times}
\newcommand{\otm}{\otimes}

\newcommand{\GAP}{{\sf GAP}}

\renewcommand{\leq}{\leqslant}
\renewcommand{\geq}{\geqslant}
\renewcommand{\spmid}{\mid}
\renewcommand{\spnmid}{\nmid}

\begin{document}
\raggedbottom
\pagestyle{myheadings}
\markboth{}{}
\thispagestyle{empty}
\setcounter{page}{1}
\begin{center} \Large\bf
On tensor products of path algebras of type $A$ \vspace*{1em} \\
\large\rm
Lutz Hille and J\"urgen M\"uller \vspace*{1em} \\
\end{center}

\begin{abstract} \noindent
We derive a formula for the Coxeter polynomial of the $s$-fold 
tensor product $\bigotimes_{i=1}^s F[\overrightarrow{A}_{n_i-1}]$
of path algebras of linearly oriented quivers of Dynkin type $A_{n_i-1}$, 
in terms of the weights $n_1\ld n_s\geq 2$, and show that conversely 
the weights can be recovered from the Coxeter polynomial of the 
tensor product.
Our results have applications in singularity theory,
in particular these algebras occur as endomorphism algebras
of tilting objects in certain stable categories of vector bundles.
\end{abstract}

\section{Introduction and results}\label{intro}

\abs
This paper is motivated by recent work \cite{KusLenMelI,KusLenMelII} 
on weighted projective lines and triangle singularities:
Let $\bX=\bX(n_1,n_2,n_3)$ be the weighted projective line,
over an algebraically closed field $F$, with respect to the
weight triple $n_1,n_2,n_3\geq 2$. Associated with $\bX$ 
there is the stable category $\un{\text{vect}}\text{-}\bX$ of 
vector bundles on $\bX$, which 
by \cite[Thm.B, Thm.6.1]{KusLenMelII} turns out to have a tilting object 
whose endomorphism algebra is isomorphic to the tensor product of algebras
$$   F[\overrightarrow{A}_{n_1-1}] \otm F[\overrightarrow{A}_{n_2-1}] 
\otm F[\overrightarrow{A}_{n_3-1}] ,$$
where $F[\overrightarrow{A}_{n_i-1}]$ denotes the path algebra
of the linearly oriented quiver of Dynkin type $A_{n_i-1}$;
for the particularly important special case of weight triples $[2,3,n]$ 
see also \cite[Prop.5.5]{KusLenMelI}.
In view of \cite[Rem.5.10, Sect.7.3]{KusLenMelII}, 
the Coxeter transformations
of these tensor product algebras are of particular interest.

\abs
It is tempting to ask whether there is a similar
interpretation of the $s$-fold tensor product 
$\bigotimes_{i=1}^s F[\overrightarrow{A}_{n_i-1}]$,
for weights $n_1\ld n_s\geq 2$ and arbitrary $s\in\N$. Indeed, 
the following interpretation of $\un{\text{vect}}\text{-}\bX$ 
as a singularity category generalises:
We first consider the hypersurface singularity
$$ S_{[n_1,n_2,n_3]}:=
   \frac{F[X_1,X_2,X_3]}{\lr{X_1^{n_1}+X_2^{n_2}+X_3^{n_3}}} ,$$
which is graded with respect to the abelian group
$\frac{\Z^4}{\lr{[n_1,0,0,1],[0,n_2,0,1],[0,0,n_3,1]}_\Z}$.
Letting $D^b(\modcat{\text{gr}}{S_{[n_1,n_2,n_3]}})$ 
be the bounded derived category of the
category of finitely generated graded $S_{[n_1,n_2,n_3]}$-modules, and 
$D^b(\projcat{\text{gr}}{S_{[n_1,n_2,n_3]}})$ 
be its full triangulated subcategory of perfect complexes, 
the associated singularity category is defined as the quotient 
$D^b_{\text{sg}}(S_{[n_1,n_2,n_3]}):=
   \frac{D^b(\modcat{\text{gr}}{S_{[n_1,n_2,n_3]}})}
        {D^b(\projcat{\text{gr}}{S_{[n_1,n_2,n_3]}})}$.
Then, by \cite{GeigleLenzing} and \cite{Orlov},
the category $\un{\text{vect}}\text{-}\bX$ is equivalent,
as a triangulated category, to the singularity category 
$D^b_{\text{sg}}(S_{[n_1,n_2,n_3]})$;
see also \cite[App.C]{KusLenMelII}.

\abs
Now, given weights $n_1\ld n_s\geq 2$, the Brieskorn-Pham
hypersurface singularity
$$ S_{[n_1\ld n_s]}:=  
     \frac{F[X_1\ld X_s]}{\lr{X_1^{n_1}+\cdots+X_s^{n_s}}} ,$$
which again is naturally graded, similarly gives rise to 
an associated singularity category
$D^b_{\text{sg}}(S_{[n_1\ld n_s]}):=
   \frac{D^b(\modcat{\text{gr}}{S_{[n_1\ld n_s]}})} 
        {D^b(\projcat{\text{gr}}{S_{[n_1\ld n_s]}})}$.
Then, by \cite[Thm.1.2]{FutakiUeda}, we have an 
equivalence of triangulated categories
$$ D^b_{\text{sg}}(S_{[n_1\ld n_s]}) \cong
   D^b\left(\modcat{}
            {\bigotimes_{i=1}^s F[\overrightarrow{A}_{n_i-1}]}\right) .$$

\abs
Finally, we mention that tensor products
$F[\overrightarrow{A}_{n_1-1}]\otm F[\overrightarrow{A}_{n_2-1}]$
have been considered, from a representation theoretic perspective, 
in \cite{Ladkani}. Thus, having all of this in mind, 
in the present paper we are interested in the algebras
$\bigotimes_{i=1}^s F[\overrightarrow{A}_{n_i-1}]$,
where $n_1\ld n_s\geq 2$ and $s\in\N$,
in their bounded derived module categories,
and in particular in their Coxeter transformations.

\AbsTT{Coxeter transformations.}
Since our investigations will take place in the realm of
finite-dimensional algebras no deep knowledge about singularity
theory or weighted projective lines will be required in the sequel. 
But we assume the reader familiar with the notion of
Coxeter transformations and their relevance in representation theory;
we are content with briefly recalling the setting:

\abs
From now on let $F$ be an arbitrary field, and let 
$A$ be a finite-dimensional (associative unital) $F$-algebra
of finite global dimension.
Let $\modcat{}{A}$ be the category of finite-dimensional (right) 
$A$-modules, let $D^b(\modcat{}{A})$ be the associated 
bounded derived category. Then their Grothendieck groups
$K_0(\modcat{}{A})$ and $K_0(D^b(\modcat{}{A}))$, respectively,
become bilinear $\Z$-lattices with respect to the Euler form,
and as such can be naturally identified.
Let $\tau_A$ be the Auslander-Reiten translation 
associated with $A$, which is a derived auto-equivalence of $A$, 
that is an auto-equivalence
of $D^b(\modcat{}{A})$ as a triangulated category.
Hence $\tau_A$ induces an isomorphism $\Phi_A$ of $K_0(\modcat{}{A})$,
being called the Coxeter transformation of $A$. 
Writing $\Phi_A$ as an element of $\GL_l(\Z)$, where
$l\in\N$ is the number of simple $A$-modules, gives rise to the 
characteristic polynomial $\chi_A\in\Z[X]$ of $\Phi_A$, 
being called the Coxeter polynomial of $A$.
Hence in particular the equivalence class of $\Phi_A$, as well as
$\chi_A$ are derived invariants of $A$, that is invariants of 
$D^b(\modcat{}{A})$ under equivalences of triangulated categories.

\abs
Of particular interest are spectral properties of $\Phi_A$:
For example, the fixed space of $\Phi_A$ coincides with the radical 
of the Euler form; and the trace of $\Phi_A$ is related to the 
Hochschild cohomology of $A$. Moreover, by Kronecker's Theorem, 
$\Phi_A$ is periodic if and only if $\Phi_A$ is diagonalisable and 
has spectral radius at most $1$. In particular, if $A$ is fractionally
Calabi-Yau, then the order of periodicity of $\Phi_A$
is related to the Calabi-Yau dimension $\frac{a}{b}$ of $A$ as follows:
Recalling that the Serre functor of $D^b(\modcat{}{A})$
induces the isomorphism $-\Phi_A$ of $K_0(\modcat{}{A})$, we get
$(\Phi_A)^b=(-\id)^{a+b}$, implying that the order of periodicity of 
$\Phi_A$ divides $2b$.

\abs
Still, despite their importance, Coxeter transformations and 
Coxeter polynomials are not too well understood, where for more 
details we refer the reader, for example, to
\cite{Happel,Lenzing,LenzingPenaI,LenzingPenaII}.
 
\AbsTT{Kronecker products.}
Given $s\in\N$ and $n_1\ld n_s\geq 2$, we now consider the
tensor product $A:=\bigotimes_{i=1}^s F[\overrightarrow{A}_{n_i-1}]$.
This of course, up to 
isomorphism, only depends on the multiset $[n_1\ld n_s]$,
but not on the order of the tensor factors; moreover, it
is independent from adding or deleting tensor
factors $F[\overrightarrow{A}_1]$, hence we will additionally  
assume that $n_i\geq 3$, for all $\{1\ld s\}$, whenever appropriate.

\abs
Anyway, for $n\geq 2$ let $\Phi_{[n]}\in\GL_{n-1}(\Z)$
be the Coxeter transformation of the path algebra 
$F[\overrightarrow{A}_{n-1}]$, 
and let $\chi_{[n]}\in\Z[X]$ be its Coxeter polynomial.
Letting 
$$ \Phi_{[n_1\ld n_s]}:=\bigotimes_{i=1}^s\Phi_{[n_i]}\in\GL_l(\Z) ,$$
where $l:=\prod_{i=1}^s(n_i-1)$,
the Coxeter transformation $\Phi_A$ associated with the tensor product
$A=\bigotimes_{i=1}^s F[\overrightarrow{A}_{n_i-1}]$ is given as 
$\Phi_A=(-1)^{s-1}\cdot\Phi_{[n_1\ld n_s]}\in\GL_l(\Z)$.
Thus we consider the map $\Phi_{[n_1\ld n_s]}$ first,
and subsequently take care of signs in order to
come back to the Coxeter transformation $\Phi_A$:
  
\abs
Let $\chi_{[n_1\ld n_s]}\in\Z[X]$ be the characteristic polynomial of 
$\Phi_{[n_1\ld n_s]}$. From $\Phi_{[2]}=[-1]\in\GL_1(\Z)$ we get
$\Phi_{[n_1\ld n_s]}=-\Phi_{[n_1\ld n_s,2]}=\Phi_{[n_1\ld n_s,2,2]}$;
as an aside, in view of the interpretation in terms of
Brieskorn-Pham hypersurface singularities, we mention that 
this observation is related to Kn\"orrer periodicity, see \cite{Knoerrer}.
Anyway, the above equalities yield
$$ \chi_{[n_1\ld n_s,2]}=(-1)^l\cdot\chi_{[n_1\ld n_s]}(-X)
   \quad\text{ and }\quad\chi_{[n_1\ld n_s,2,2]}=\chi_{[n_1\ld n_s]} .$$ 
Thus we may additionally assume that $n_i\geq 3$, for all $\{1\ld s-1\}$, 
that is the multiset $[n_1\ld n_s]$ contains the element $2$ at most once,
whenever appropriate.

\abs
Our main results now are a formula expressing  
the polynomial $\chi_{[n_1\ld n_s]}\in\Z[X]$ 
as a rational function in $\Q(X)$, and a recognition result:

\Thm\label{thm1}
Let $[n_1\ld n_s]$ be a multiset, and for any subset 
$\emp\neq\cJ\sseq\cI:=\{1\ld s\}$ let 
$n_\cJ:=\lcm(n_i;i\in\cJ)\in\N$, and let $n_\emp:=1$. 
Then we have
$$ \chi_{[n_1\ld n_s]}=\prod_{\cJ\sseq\cI}\left(
   \left(X^{n_\cJ}-1\right)^{\frac{\prod_{i\in\cJ}n_i}{n_\cJ}}
   \right)^{(-1)^{s-|\cJ|}} .$$ 

\Thm\label{thm2}
Let $[n_1\ld n_s]$ and $[\wti{n}_1\ld\wti{n}_{\wti{s}}]$ be multisets
each containing the element $2$ at most once. Then we have 
$\chi_{[n_1\ld n_s]}=\chi_{[\wti{n}_1\ld\wti{n}_{\wti{s}}]}\in\Z[X]$
if and only if the multisets 
$[n_1\ld n_s]$ and $[\wti{n}_1\ld\wti{n}_{\wti{s}}]$ coincide.

\abs\abs
Actually, in order to prove Theorem \ref{thm2} we prove much more: 
We give an algorithm to recover the multiset $[n_1\ld n_s]$,
provided it contains the element $2$ at most once,
from the characteristic polynomial $\chi_{[n_1\ld n_s]}$ alone.
Both the algorithm, and the formula in Theorem \ref{thm1} 
to compute $\chi_{[n_1\ld n_s]}$ from $[n_1\ld n_s]$,
are easily implemented, for example into the 
computer algebra system \GAP{} \cite{GAP}; this has been used
to verify the examples given in \ref{ex1}, \ref{ex1b}, 
\ref{ex2} and \ref{ex3}.

\abs
Having these general conclusions, and the techniques to prove them,
at hand, we are also able to derive a few further properties of 
the spectrum of $\Phi_{[n_1\ld n_s]}$, which by Theorem \ref{thm1} 
consists of roots of unity of order dividing $n_\cI$.
Firstly, we deal with the question whether the roots of unity
of highest possible order actually are elements of the spectrum:
 
\Cor\label{cor1}
Let $[n_1\ld n_s]$ be a multiset containing the element $2$ at most once.
Then the order of periodicity of $\Phi_{[n_1\ld n_s]}$ equals 
$n_\cI$, where moreover all primitive $n_\cI$-th roots of unity
are eigenvalues of $\Phi_{[n_1\ld n_s]}$.

\abs\abs
Moreover, we are able to give a general criterion to decide when $1$,
that is the root of unity of lowest possible order,
is an element of the spectrum; due to its technicality
it is only stated in Theorem \ref{m1thm}. Finally, a consideration
of the multiplicity of $1$ as an element of the spectrum leads
to the following result: Recall that for any $f\in\Z[X]\smin\{0\}$ 
the reciprocal polynomial is defined as 
$f^\ast:=X^{\deg(f)}\cdot f(X^{-1})\in\Z[X]\smin\{0\}$,
and $f$ is called self-reciprocal if $f^\ast=f$.

\Cor\label{cor2}
Let $[n_1\ld n_s]$ be a multiset. Then 
$\chi_{[n_1\ld n_s]}^\ast=\pm\chi_{[n_1\ld n_s]}$,
where $\chi_{[n_1\ld n_s]}$ is self-reciprocal if and only if
$\gcd(s,n_1\ld n_s)$ is odd.

\AbsTT{Back to tensor products.}
Returning to the original motivation, 
we draw a few immediate conclusions about the tensor product 
$A:=\bigotimes_{i=1}^s F[\overrightarrow{A}_{n_i-1}]$,
where we now assume that $n_i\geq 3$ for all $i\in\{1\ld s\}$:

\abs\bfa
The associated Coxeter polynomial being given as
$$ \left\{\begin{array}{ll}
\chi_{[n_1\ld n_s]},&\text{ if $s$ is odd}, \\
\chi_{[n_1\ld n_s,2]},&\text{ if $s$ is even}, \\
\end{array}\right. $$
by Corollary \ref{cor1} the Coxeter transformation $\Phi_A$
is periodic of order
$$ \left\{\begin{array}{rclll}
n_\cI&=&\lcm(n_1\ld n_s),&\text{ if $s$ is odd}, \\
\lcm(n_\cI,2)&=&\lcm(n_1\ld n_s,2),&\text{ if $s$ is even}. \\
\end{array}\right. $$
In particular, for $s=3$ we recover \cite[Prop.7.6]{KusLenMelII};
note that the proof given there heavily depends on using
weighted projective lines.

\abs\bfb
Moreover, since Coxeter polynomials are derived invariants,
from Theorem \ref{thm2} we infer that tensor products 
$A=\bigotimes_{i=1}^s F[\overrightarrow{A}_{n_i-1}]$ and 
$\bigotimes_{i=1}^{\wti{s}} F[\overrightarrow{A}_{\wti{n}_i-1}]$,
where $\wti{s}\in\N$, and $\wti{n}_i\geq 3$ for all $i\in\{1\ld\wti{s}\}$, 
are derived equivalent if and only if the multisets
$[n_1\ld n_s]$ and $[\wti{n}_1\ld\wti{n}_{\wti{s}}]$ coincide.

\abs\bfc
Finally, we note that the algebra
$A=\bigotimes_{i=1}^s F[\overrightarrow{A}_{n_i-1}]$
is fractionally Calabi-Yau, see \cite[Sect.1.4]{Ladkani},
where since $F[\overrightarrow{A}_{n_i-1}]$ has Calabi-Yau dimension
$\frac{n_i-2}{n_i}$, see \cite[Thm.4.1]{MiyachiYekutieli},
we infer that $A=\bigotimes_{i=1}^s F[\overrightarrow{A}_{n_i-1}]$ 
has Calabi-Yau dimension
$$ \frac{\sum_{i=1}^s \frac{n_\cI}{n_i}\cdot(n_i-2)}{n_\cI}
  =\frac{s\cdot n_\cI-2\cdot\sum_{i=1}^s \frac{n_\cI}{n_i}}{n_\cI} ;$$
for the case $s=3$ see also \cite[Prop.7.5]{KusLenMelII}.

\abs
This yields $(\Phi_A)^{n_\cI}=(-\id)^{(s-1)\cdot n_\cI}$,
just implying the obvious fact that the order of periodicity of $\Phi_A$ 
divides $n_\cI$ respectively $\lcm(n_\cI,2)$, 
whenever $s$ is odd respectively even.
But comparing with the actual order of periodicity of $\Phi_A$ as determined
above, we observe that we indeed detect the fractional part
of the Calabi-Yau dimension of 
$A=\bigotimes_{i=1}^s F[\overrightarrow{A}_{n_i-1}]$
in terms of the associated Coxeter transformation.

\AbsTT{Outline.}
This paper is organised as follows:
In Section \ref{coxpol} we prove Theorem \ref{thm1};
in order to do so, we apply character theory of finite abelian groups,
where we recall the necessary facts, but for more details refer 
the reader, for example, to \cite{Isaacs}.
In Section \ref{findfactors} we prove Theorem \ref{thm2}
and the first half of Corollary \ref{cor1};
in order to do so, we are going to apply lattice theoretic
M\"obius inversion, where again we recall the necessary facts,
but for more details refer the reader, for example, to \cite{Rota}.
In Section \ref{coxspec} we finally prove the second half of 
Corollary \ref{cor1}, state and prove Theorem \ref{m1thm}, and prove
Corollary \ref{cor2}.

\AbsTT{Acknowledgement.}
We thank Helmut Lenzing for drawing our attention to these
questions, for various enlightening discussions, and in 
particular for sharing his thoughts about the relevance 
of tensor product algebras in singularity theory with us.
We greatfully acknowledge financial support by the
German Research Council (DFG) in the framework of the 
Scientific Priority Program (SPP 1388) `Representation Theory'.
The first author thanks the Max Planck Institute for 
Mathematics (MPI) in Bonn for its hospitality.

\section{Characteristic polynomials of Coxeter transformations}\label{coxpol}

\abs
We proceed towards a proof of Theorem \ref{thm1}, where we from now on 
again allow for $n_i\geq 2$ for all $i\in\{1\ld s\}$.

\Abs
We recall the determination of the Coxeter polynomial of the
path algebra $F[\overrightarrow{A}_{n-1}]$, which is of course
well-known, see \cite[Ch.VIII.5]{AusReiSma}: Let 
$$ \cC_{[n]}:=\begin{bmatrix}
1&1&1&\ldots&1\\
&1&1&\ldots&1\\
&&1&\ldots&1\\
&&&\ddots&\vdots\\
&&&&1\\
\end{bmatrix}\in\GL_{n-1}(\Z) $$
be the Cartan matrix associated with $F[\overrightarrow{A}_{n-1}]$,
describing the classes of the projective indecomposable modules 
in $K_0(F[\overrightarrow{A}_{n-1}])\cong\Z^{n-1}$ with respect
to the `standard' basis consisting of the classes of the simple modules;
the injective indecomposable modules are similarly 
described by $\cC_{[n]}^{\tr}$. 
Hence, by \cite[Ch.VIII.2]{AusReiSma}, the Coxeter transformation
$\Phi_{[n]}\in\Z^{(n-1)\tm(n-1)}$ is uniquely determined by the equation
$\cC_{[n]}\cdot\Phi_{[n]}=-\cC_{[n]}^{\tr}$, thus
$$ \Phi_{[n]}=-\cC_{[n]}^{-1}\cdot\cC_{[n]}^{\tr}
=\begin{bmatrix}
.&1&&&\\
&.&1&&\\
&&\ddots&\ddots\\
&&&.&1\\
-1&-1&\ldots&-1&-1\\
\end{bmatrix}\in\GL_{n-1}(\Z) ;$$
in particular, we have $\Phi_{[2]}=[-1]\in\GL_1(\Z)$.

\abs
Since $\Phi_{[n]}$ is described by a companion matrix, 
its characteristic polynomial $\chi_{[n]}\in\Z[X]$ is given as
$$ \chi_{[n]}
=\sum_{i=0}^{n-1}X^i=\frac{X^n-1}{X-1}
=\sum_{i=1}^{n-1}(X-\zt_n^i)\in\C[X] ,$$
where 
$\zt_n:=\exp(\frac{2\pi\sqrt{-1}}{n})\in\C$ 
is the standard primitive complex $n$-th root of unity.
Hence $\Phi_{[n]}$ is diagonalisable over $\C$,
is periodic of order $n$, 
and its spectrum is given by the pairwise distinct eigenvalues 
$\{\zt_n^i\in\C;i\in(\Z/n)\smin\{0\}\}$.
 
\Abs
Since the Cartan matrix of the tensor product 
$\bigotimes_{i=1}^s F[\overrightarrow{A}_{n_i-1}]$
is given by the Kronecker product $\bigotimes_{i=1}^s\cC_{[n_i]}$,
the associated Coxeter transformation is 
$$ -\left(\bigotimes_{i=1}^s\cC_{[n_i]}\right)^{-1}\cdot
\left(\bigotimes_{i=1}^s\cC_{[n_i]}\right)^{\tr}
=-\bigotimes_{i=1}^s\left(\cC_{[n_i]}^{-1}\cdot\cC_{[n_i]}^{\tr}\right)
=(-1)^{s-1}\cdot\bigotimes_{i=1}^s\Phi_{[n_i]} .$$

\abs
We let $\Phi_{[n_1\ld n_s]}:=\bigotimes_{i=1}^s\Phi_{[n_i]}$.
By the above, $\Phi_{[n_1\ld n_s]}$ is also diagonalisable over $\C$ 
and periodic, its complex eigenvalues being given as the multiset
$$\cM:=\left\{
\prod_{i=1}^s\zt_{n_i}^{a_i}\in\C;a_i\in(\Z/{n_i})\smin\{0\}\right\} .$$
In other words, its characteristic polynomial $\chi_{[n_1\ld n_s]}\in\Z[X]$ 
is given as 
$$ \chi_{[n_1\ld n_s]}
=\prod_{[a_1\ld a_s]\in\prod_{i=1}^s((\Z/{n_i})\smin\{0\})}
 \left(X-\prod_{i=1}^s\zt_{n_i}^{a_i}\right)\in\C[X] .$$
Note that these considerations are reminiscent of the
approach taken in \cite{Happel}.

\abs
Hence in particular $\Phi_{[n_1\ld n_s]}$ has spectral radius $1$,
and since the above description of the Coxeter transformation 
yields $\det((-1)^{s-1}\cdot\Phi_{[n_1\ld n_s]})=(-1)^l$,
where $l:=\prod_{i=1}^s(n_i-1)$, we for later use note that
$$ \det(\Phi_{[n_1\ld n_s]})=(-1)^{ls} .$$
We now proceed towards 
the asserted explicit formula for $\chi_{[n_1\ld n_s]}$:

\Abs
To this end, 
let $G_n:=\lr{z}$ be a (multiplicative) cyclic group of order $n\geq 2$,
and let $z$ be a generator. Moreover, let
$G^\ast_n:=\{(\zt^\ast_n)^j;j\in\Z/n\}$ be its (multiplicative)
character group, where $\zt^\ast_n\cn G_n\ra\C^\ast\cn z\mt\zt_n$,
then the integral group ring $\Z[G^\ast_n]$, consisting
of the formal $\Z$-linear combinations of the characters
in $G^\ast_n$, is called the (additive) group of generalised
characters of $G_n$. 

\abs
In particular, let $1_n:=(\zt^\ast_n)^0\in G^\ast_n$ be the 
trivial character, that is the character of the representation 
$G_n\ra\GL_1(\C)\cn z\mt 1$,
and let $\rho_n:=\sum_{j=0}^{n-1}(\zt^\ast_n)^j\in\Z[G^\ast_n]$
be the regular character, that is the character of the regular 
action of $G_n$ on the complex group algebra $\C[G_n]$.
Then we have $\rho_n(1)=n$ and 
$\rho_n(x)=0$ for all $x\in G_n\smin\{1\}$,
and conversely for any generalised character $\chi\in\Z[G^\ast_n]$ 
such that $\chi(x)=0$ for all $x\in G_n\smin\{1\}$
we have $\chi=\frac{\chi(1)}{n}\cdot\rho_n$.

\abs
Now, by mapping $z\mt\Phi_{[n]}$ we get a faithful representation 
of $G_n\ra\GL_{n-1}(\C)$, whose character is given as
$\ph_n:=\sum_{j=1}^{n-1}(\zt^\ast_n)^j=\rho_n-1_n\in\Z[G^\ast_n]$.

\Abs
Now we consider the direct product
$$ G:=\prod_{i=1}^s G_{n_i}=\prod_{i=1}^s\lr{z_i} $$
of cyclic groups $G_{n_i}=\lr{z_i}$ with chosen generators $z_i$.
Let $\rho\in\Z[G^\ast]$ be the regular character of $G$,
which hence can be written as the outer tensor product 
$\rho=\rho_{n_1}\otm\cdots\otm\rho_{n_s}$ 
of the regular characters $\rho_{n_i}\in\Z[G_{n_i}^\ast]$.
Moreover, let $\ph:=\ph_{n_1}\otm\cdots\otm\ph_{n_s}\in\Z[G^\ast]$ be 
the outer tensor product of the characters $\ph_{n_i}\in\Z[G_{n_i}^\ast]$. 
Hence we have 
$\ph=(\rho_{n_1}-1_{n_1})\otm\cdots\otm(\rho_{n_s}-1_{n_s})\in\Z[G^\ast]$.

\abs
Considering the representation of $G$ affording $\ph$,
since the tensor factors $\ph_{n_i}$ are afforded by the
representations $z_i\mt\Phi_{[n_i]}$,
we infer that the element $z:=\prod_{i=1}^s z_i\in G$
is represented by the matrix 
$\Phi_{[n_1\ld n_s]}=\bigotimes_{i=1}^s\Phi_{[n_i]}$.
Hence we consider the cyclic subgroup
$$ H:=\lr{z}=\lr{\prod_{i=1}^s z_i}\leq G ,$$
which in particular has order $n_\cI=\lcm(n_1\ld n_s)$.
Thus the restriction $\ph|_H$ can be viewed as an element of $\Z[H^\ast]$,
and hence may be written as a $\Z$-linear combination of linear 
characters of $H$. This will yield the characteristic polynomial 
$\chi_{[n_1\ld n_s]}$ of $\Phi_{[n_1\ld n_s]}$ 
as a rational function in $\C(X)$,
where the linear characters occurring, evaluated at $z\in H$,
describe its zeroes and poles.

\Abs
To this end, for $\cJ\sseq\cI$ let 
$$ G_\cJ:=\prod_{i\in\cJ}G_{n_i}=\prod_{i\in\cJ}\lr{z_i} \leq G $$
be the subgroup generated by the direct factors indicated by $\cJ$.
Hence we have a group epimorphism
$$ \al_\cJ\cn G\ra G_\cJ\cn\left\{\begin{array}{ll}
z_i\mt z_i, &\text{ if }i\in\cJ, \\
z_i\mt 1, &\text{ if }i\not\in\cJ. \\
\end{array}\right. $$
Let $\rho_\cJ=\bigotimes_{i\in\cJ}\rho_{n_i}\in\Z[G_\cJ^\ast]$
be the regular character of $G_\cJ$, 
and let $\wh{\rho}_\cJ:=\rho_\cJ\circ\al_\cJ\in\Z[G^\ast]$ 
be its inflation to $G$ via $\al_\cJ$.
Thus expanding the outer tensor product
$\ph=(\rho_{n_1}-1_{n_1})\otm\cdots\otm(\rho_{n_s}-1_{n_s})$ 
by distributivity we get
$$ \ph=\sum_{\cJ\sseq\cI}(-1)^{s-|\cJ|}\cdot(\rho_\cJ\circ\al_\cJ) 
=\sum_{\cJ\sseq\cI}(-1)^{s-|\cJ|}\cdot\wh{\rho}_\cJ \in\Z[G^\ast] .$$

\Abs
Finally, let $z_\cJ:=\prod_{i\in\cJ} z_i$ and  
$$ H_\cJ:=\lr{z_\cJ}\leq G_\cJ ,$$
which in particular has order $n_\cJ=\lcm(n_i;i\in\cJ)$.
Then $\al_\cJ$ restricts to an epimorphism 
$H\ra H_\cJ\cn z\mt z_\cJ$, which we denote by the same symbol.
Let $\rho_{n_\cJ}\in\Z[H_\cJ^\ast]$ be the regular character of $H_\cJ$, 
and let $\wh{\rho}_{n_\cJ}:=\rho_{n_\cJ}\circ\al_\cJ\in\Z[H^\ast]$ 
be its inflation to $H$ via $\al_\cJ$. Then letting 
$\zt^\ast_{n_\cJ}:=(\zt^\ast_{n_\cI})^{\frac{n_\cI}{n_\cJ}}\in\Z[H^\ast]$
we have
$$ \wh{\rho}_{n_\cJ}=\sum_{j=0}^{n_\cJ-1}(\zt^\ast_{n_\cJ})^j\in\Z[H^\ast] .$$

\abs
Now we have $\rho_\cJ(1)=|G_\cJ|=\prod_{i\in\cJ}n_i$
and $\rho_{n_\cJ}(1)=|H_\cJ|=n_\cJ$.
Since $\rho_\cJ(x)=0$ for all $x\in G_\cJ\smin\{1\}$, we conclude that
$$ \rho_\cJ|_{H_\cJ}
=\frac{\prod_{i\in\cJ}n_i}{n_\cJ}\cdot\rho_{n_\cJ}\in\Z[H_\cJ^\ast] .$$
Hence for all $x\in H$ we have
$$ \wh{\rho}_\cJ(x)=\rho_\cJ(\al_\cJ(x))
=\frac{\prod_{i\in\cJ}n_i}{n_\cJ}\cdot\rho_{n_\cJ}(\al_\cJ(x))
=\frac{\prod_{i\in\cJ}n_i}{n_\cJ}\cdot\wh{\rho}_{n_\cJ}(x) ,$$
that is
$$ \wh{\rho}_\cJ|_H
=\frac{\prod_{i\in\cJ}n_i}{n_\cJ}\cdot\wh{\rho}_{n_\cJ}\in\Z[H^\ast] .$$
This implies
$$ \ph|_H=\sum_{\cJ\sseq\cI}(-1)^{s-|\cJ|}\cdot\wh{\rho}_\cJ|_H 
=\sum_{\cJ\sseq\cI}\left((-1)^{s-|\cJ|}\cdot
 \frac{\prod_{i\in\cJ}n_i}{n_\cJ}\cdot\wh{\rho}_{n_\cJ}\right)
\in\Z[H^\ast] ,$$
where in turn
$\wh{\rho}_{n_\cJ}=\sum_{j=0}^{n_\cJ-1}(\zt^\ast_{n_\cJ})^j$.
Hence evaluating at $z\in H$ we get
$$ \chi_{[n_1\ld n_s]}=\prod_{\cJ\sseq\cI}\left( \left(
   \prod_{j=0}^{n_\cJ-1}(X-\zt^j_{n_\cJ})
   \right)^{(-1)^{s-|\cJ|}\cdot\frac{\prod_{i\in\cJ}n_i}{n_\cJ}}\right)
   \in\C(X) ,$$ 
which is the desired description, proving Theorem \ref{thm1}.
\QED

\Ex\label{ex1}
We write out the above formula explicitly for small values of $s\in\N$,
and verify that the polynomials obtained indeed behave as expected:

\abs\bfa
For $s=1$, letting $n:=n_1$, we of course just recover 
$\chi_{[n]}=\frac{X^n-1}{X-1}$.

\abs\bfb
For $s=2$, letting $n_{12}:=\lcm(n_1,n_2)$ and
$g_{12}:=\gcd(n_1,n_2)=\frac{n_1 n_2}{n_{12}}$, we have
$$ \chi_{[n_1,n_2]}
=\frac{(X-1)(X^{n_{12}}-1)^{g_{12}}}
      {(X^{n_1}-1)(X^{n_2}-1)} ;$$
in particular, if $n_1$ and $n_2$ are coprime this yields
$$ \chi_{[n_1,n_2]}=\frac{(X-1)(X^{n_1n_2}-1)}
      {(X^{n_1}-1)(X^{n_2}-1)} .$$

\abs
Even more specifically, for the multiset $[n,2]$ we get 
$\chi_{[n,2]}=\frac{X^n-(-1)^n}{X+1}$, 
thus in any case we indeed have 
$\chi_{[n,2]}=(-1)^{n-1}\cdot\chi_{[n]}(-X)$ as expected.

\abs\bfc
For $s=3$, abbreviating $n_{ij}:=\lcm(n_i,n_j)$ and
$g_{ij}:=\gcd(n_i,n_j)=\frac{n_i n_j}{n_{ij}}$,
for $i<j\in\{1,2,3\}$, and $g_{123}:=\gcd(n_1,n_2,n_3)$,
and noting that the principle of inclusion-exclusion yields
$n_{123}:=\lcm(n_1,n_2,n_3)=\frac{n_1 n_2 n_3 g_{123}}{g_{12}g_{13}g_{23}}$,
we have
$$ \chi_{[n_1,n_2,n_3]}=
\frac{(X^{n_1}-1)(X^{n_2}-1)(X^{n_3}-1)
      (X^{n_{123}}-1)^{\frac{g_{12}g_{13}g_{23}}{g_{123}}}}
      {(X-1)(X^{n_{12}}-1)^{g_{12}}(X^{n_{13}}-1)^{g_{13}}
       (X^{n_{23}}-1)^{g_{23}}} ;$$
in particular, if $n_1$, $n_2$ and $n_3$ are pairwise coprime this yields
$$ \chi_{[n_1,n_2,n_3]}=
\frac{(X^{n_1}-1)(X^{n_2}-1)(X^{n_3}-1)(X^{n_1n_2n_3}-1)}
      {(X-1)(X^{n_1n_2}-1)(X^{n_1n_3}-1)(X^{n_2n_3}-1)} .$$

\abs
Even more specifically, for the multiset $[n,2,2]$ we get 
$\chi_{[n,2,2]}=\frac{X^n-1}{X-1}$, independently
of whether $n$ is odd or even;
thus we indeed have $\chi_{[n,2,2]}=\chi_{[n]}$ as expected.
Moreover, for $[n_1,n_2,2]$, by a similar case distinction
with respect to the parity of $n_1$ and $n_2$, we indeed get, as expected,
$$ \chi_{[n_1,n_2,2]}
   =(-1)^{(n_1-1)(n_2-1)}\cdot\chi_{[n_1,n_2]}(-X) .$$

\Ex\label{ex1b}
We now consider a few examples which play a particular role in the
context of weighted projective lines:

\abs\bfa
For the tensor product 
$F[\overrightarrow{A}_2]\otm F[\overrightarrow{A}_{n-1}]$, where $n\geq 2$, 
which is associated with the most prominent series  $\bX(2,3,n)$ of
weighted projective lines, see \cite{KusLenMelI},
we get the Coxeter polynomial
$$ \chi_{[2,3,n]}=\chi_{[3,n]}(-X)=\left\{\begin{array}{ll}
\frac{(X+1)(X^{3n}-(-1)^{n})}{(X^3+1)(X^n-(-1)^{n})},
   &\text{ if }3\spnmid n, \\
\frac{(X+1)(X^n-(-1)^{n})^2}{(X^3+1)},
   &\text{ if }3\spmid n. \rule{0em}{1.5em} \\
\end{array}\right. $$

\abs
In particular, for $n\in\{2\ld 5\}$ we get 
$$ \chi_{[2,3,n]}=\left\{\begin{array}{ll}
\Phi_3,&\text{ if }n=2, \\
\Phi_2^2\Phi_6,&\text{ if }n=3, \rule{0em}{1.5em} \\
\Phi_3\Phi_{12},&\text{ if }n=4, \rule{0em}{1.5em} \\
\Phi_{30},&\text{ if }n=5, \rule{0em}{1.5em} \\
\end{array}\right. $$
where $\Phi_d\in\Z[X]$, for $d\in\N$,
denotes the $d$-th cyclotomic polynomial.
These coincide with the Coxeter polynomials 
of the path algebras of the quivers of Dynkin type
$A_2$, $D_4$, $E_6$, and $E_8$, respectively, see for example 
\cite[Sect.18.5]{Lenzing}.

\abs
This was to be expected, since, apart from $\bX(2,2,n)$,
the above cases are precisely those of
weighted projective lines having positive Euler characteristic, 
see \cite[Sect.5.1]{KusLenMelII}, and the associated tensor products are
by the ADE chain phenomenon derived equivalent to the path algebras 
named, see \cite[Prop.5.5, Prop.5.16]{KusLenMelI},
or independently \cite[Cor.1.2]{Ladkani}.

\abs\bfb
Moreover, for the tensor products
$$ F[\overrightarrow{A}_2]\otm F[\overrightarrow{A}_2]
                          \otm F[\overrightarrow{A}_2],
\quad
F[\overrightarrow{A}_3]\otm F[\overrightarrow{A}_3],
\quad
F[\overrightarrow{A}_2]\otm F[\overrightarrow{A}_5], $$
being associated with the weighted projective lines 
$\bX(3,3,3)$, $\bX(2,4,4)$ and $\bX(2,3,6)$, respectively,
we get the Coxeter polynomials
$$ \begin{array}{lclcl}
\chi_{[3,3,3]}&=&(X-1)^2\cdot\left(\frac{X^3-1}{X-1}\right)^3
              &=&\Phi_1^2\Phi_3^3, \\
\chi_{[2,4,4]}&=&(X-1)^2\cdot\left(\frac{X^2-1}{X-1}\right)\cdot
                             \left(\frac{X^4-1}{X-1}\right)^2
              &=&\Phi_1^2\Phi_2^3\Phi_4^2, \rule{0em}{1.5em} \\
\chi_{[2,3,6]}&=&(X-1)^2\cdot\left(\frac{X^2-1}{X-1}\right)\cdot
                             \left(\frac{X^3-1}{X-1}\right)\cdot              
                             \left(\frac{X^6-1}{X-1}\right)  
&=&\Phi_1^2\Phi_2^2\Phi_3^2\Phi_6, \rule{0em}{1.5em} \\
\end{array} $$
which are of canonical type, see \cite[Prop.18.6]{Lenzing}.

\abs
This was to be expected, since, apart from $\bX(2,2,2,2)$,
the above cases are precisely those of weighted projective 
lines having zero Euler characteristic. In other words,
these are the tubular weight types,
see \cite[Sect.5.2]{KusLenMelII}, hence the associated tensor products
are derived equivalent to the canonical algebras
of the respective type, see \cite[Thm.5.6]{KusLenMelII},
or \cite[Sect.18.6.2]{Lenzing} for the case $[3,3,3]$.

\section{Recovering the tensor factors}\label{findfactors}

\abs
We now conversely assume only the characteristic polynomial
$\chi_\Phi\in\Z[X]$, where $\Phi:=\Phi_{[n_1\ld n_s]}$ is as above, 
to be known. Thus we hide the multiset from the notation,
and we aim to recover the number $s\in\N$ of tensor factors
and the multiset $[n_1\ld n_s]$ from $\chi_\Phi$ alone.
This will in particular entail Theorem \ref{thm2}, 
as soon as we additionally assume that the multiset $[n_1\ld n_s]$
contains the element $2$ at most once.

\Abs
Recall that for any $d\in\N$ the multiplicity
$m_d\in\N_0$ of the primitive $d$-th root of unity $\zt_d\in\C$
as an element of the multiset $\cM$ of complex eigenvalues of $\Phi$
is given as the order of $\zt_d$ as a zero of
$\chi_\Phi\in\C[X]$;
note that $m_d$ only depends on the order $d$, but not on the particular
choice of a primitive $d$-th root of unity.
Hence, by evaluating $\chi_\Phi$ at the standard primitive 
$d$-th root of unity $\zt_d\in\C$ for increasing values of $d\in\N$,
the multiplicities $m_d\in\N_0$, the number $\max(d\in\N;m_d>0)\in\N$,
and thus the number
$$ N:=\lcm(d\in\N;m_d>0)\in\N $$
can be determined from $\chi_\Phi$; note that $N$ 
is the order of periodicity of $\Phi$.

\abs
To facilitate the subsequent analysis, assume that $n\in\N$ 
is a multiple of the currently unknown number $n_\cI=\lcm(n_1\ld n_s)$.
Then we may rewrite $\cM$ as 
$$ \cM=\left\{\prod_{i=1}^s\zt_n^{\frac{n}{n_i}\cdot a_i}\in\C;
              a_i\in(\Z/{n_i})\smin\{0\}\right\} .$$
We have $m_d>0$ only if $d\spmid n_\cI$, hence 
we get $N\spmid n_\cI$ as well, 
but a priorly we might have $m_{n_\cI}=0$ or even $N<n_\cI$.

\abs
Still, a suitable, but crudely large $n\in\N$ can be found as follows:
Since $\deg(\chi_\Phi)=\prod_{i=1}^s(n_i-1)$
we have $n_i\leq\deg(\chi_\Phi)+1$, for all $i\in\{1\ld s\}$,
hence we may choose $n\in\N$ as a multiple of 
$\lcm(1,2\ld\deg(\chi_\Phi),\deg(\chi_\Phi)+1)$.
This a priori choice admittedly in general is much larger than $n_\cI$,
but is sufficient for the moment, 
and we will come back to that point in \ref{clarify} again.

\Abs
Converting to an additive setting, for any $d\spmid n$ 
we have the natural embedding 
$\ia^d_n\cn\Z/d\ra\Z/n\cn x\pmod{d}\mt\frac{n}{d}\cdot x\pmod{n}$.
Hence we have the natural homomorphism 
$$ \ph:=[\ia^{n_1}_n\ld\ia^{n_s}_n]\cn
   \Z/{n_1}\oplus\cdots\oplus\Z/{n_s}\ra\Z/n\cn 
   [x_1\ld x_s]\mt\sum_{i=1}^s\frac{n}{n_i}\cdot x_i .$$
Letting 
$$ \cO:=((\Z/{n_1})\smin\{0\})\oplus\cdots\oplus((\Z/{n_s})\smin\{0\})
\sseq\Z/{n_1}\oplus\cdots\oplus\Z/{n_s} ,$$
the multiset $\cM$ is described by the map $\ph|_{\cO}$, 
that is, if $x\in\Z/n$ has additive order $d\spmid n$,
then we have $|(\ph|_{\cO})^{-1}(x)|=m_d$.
Thus the task is, given this information 
on the fibres of $\ph|_{\cO}$, to recover the multiset $[n_1\ld n_s]$.

\Abs
In order to do so, we are going to apply lattice theoretic 
M\"obius inversion. We recall the necessary facts:
Let $\cL$ be a finite combinatorial lattice, that is a non-empty 
finite set equipped with a partial order $\leq$, such that for 
any $x,y\in\cL$ there is a smallest upper bound $x+y\in\cL$
and a largest lower bound $x\cap y\in\cL$ with respect to $\leq$.

\abs
Given a map $f\cn\cL\ra\C$, the associated sum function is defined as
$$ f^+\cn\cL\ra\C\cn x\mt\sum_{y\leq x}f(y) .$$ 
To recover $f$ from the knowledge of $f^+$, we may 
use the M\"obius function associated with $\cL$, 
which is recursively defined as 
$$\mu_\cL\cn\cL\tm\cL\ra\C\cn[x,y]\mt
\left\{\begin{array}{rl}
0,&\text{ if }x\not\leq y,\\
1,&\text{ if }x=y,\\
-\sum_{x\leq z<y}\mu_\cL(x,z),&\text{ if }x<y .\\
\end{array}\right. $$
Note that $\mu_\cL$ is essentially described by the identities
$\sum_{x\leq z\leq y}\mu_\cL(x,z)=0$ whenever $x<y\in\cL$. 
Then for all $x\in\cL$ we have the M\"obius inversion 
formula to recover $f$ from $f^+$, see \cite[Prop.2]{Rota},
$$ f(x)=\sum_{y\leq x}f^+(y)\mu_\cL(y,x) .$$

\abs
Moreover, restricting to an interval $\cL'=\{z\in\cL;u\leq z\leq v\}$,
where $u\leq v\in\cL$, we get the M\"obius function
$\mu_{\cL'}=\mu_\cL|_{\cL'\tm\cL'}$, see \cite[Prop.4]{Rota}; and 
letting $\cL^\ast$ be the dual lattice associated with $\cL$, obtained 
by reverting the partial order, its M\"obius function is
$\mu_{\cL^\ast}(x,y)=\mu_\cL(y,x)$, for all $x,y\in\cL$,
see \cite[Prop.3]{Rota}.

\abs
We are going to apply M\"obius inversion in two cases:
Firstly, for the partially ordered set of all subsets of a 
finite set $\cI$, where the partial order is given by
set theoretic inclusion $\sseq$, the non-zero values of 
the associated M\"obius function $\mu_\cI$ are given as
$$ \mu_\cI(\cK,\cJ)=(-1)^{|\cJ|-|\cK|}
\quad\text{ whenever }\cK\sseq\cJ\sseq\cI .$$
The resulting inversion principle is known as
the principle of inclusion-exclusion.

\abs
Secondly, given $n\in\N$, for the partially ordered set
$\{d\in\N;d\spmid n\}$ of all divisors of $n$,
where the partial order is given by the divisibility relation, 
the essential part of the associated (lattice theoretic) 
M\"obius function $\mu_n$ is given as, 
see for example \cite[Ch.1.4.9]{Bundschuh}, 
$$ \mu_n(x,y)=\mu(\frac{y}{x})\quad\text{ whenever }x\spmid y\spmid n ,$$
where $\mu\cn\N\ra\{-1,0,1\}$ is the number theoretic
M\"obius function given as
$$ \mu\cn d\mt\left\{\begin{array}{rl}
0,&\text{ if $d$ is not squarefree},\\
(-1)^t,&\text{ if $d$ is the product of $t\in\N_0$ pairwise
               distinct primes}.\\
\end{array}\right. $$

\Abs
We are now prepared for our first application of M\"obius inversion:
For $\cJ\sseq\cI:=\{1\ld s\}$ let 
$$ G_\cJ:=\bigoplus_{i\in\cJ}\Z/{n_i}\leq
   \Z/{n_1}\oplus\cdots\oplus\Z/{n_s} $$
be the subgroup generated by the direct summands indicated by $\cJ$,
and let $\ph_\cJ:=\ph|_{G_\cJ}\cn G_\cJ\ra\Z/n$.
Then, given $x\in\Z/n$, by taking the support of elements of $G_\cJ$ into
account, the cardinality $|\ph_\cJ^{-1}(x)|\in\N_0$ is given as 
$$ |\ph_\cJ^{-1}(x)|=\sum_{\cK\sseq\cJ}|(\ph|_{G_\cK\cap\cO})^{-1}(x)| .$$
Hence the principle of inclusion-exclusion, applied to the 
full set $\cI$, yields
$$ |(\ph|_{\cO})^{-1}(x)|=|(\ph|_{G_\cI\cap\cO})^{-1}(x)|
=\sum_{\cJ\sseq\cI}(-1)^{s-|\cJ|}\cdot|\ph_\cJ^{-1}(x)| ,$$
reducing the task of computing the size of fibres of $\ph|_{\cO}$
to determining the size of fibres of the homomorphisms $\ph_\cJ$,
for all $\cJ\sseq\cI$.

\abs
For $\cJ\sseq\cI$ let still $n_\cJ:=\lcm(n_i;i\in\cJ)$,
where we again set $n_\emp:=1$.
Then we have $\im(\ph_\cJ)=\im(\ia_n^{n_\cJ})\leq\Z/n$, that is
$|\im(\ph_\cJ)|=n_\cJ$. 
Moreover, if $x\in\Z/n$ has additive order $d\spmid n$, then
we have $x\in\im(\ph_\cJ)$ if and only if $d\spmid n_\cJ$.
In this case we have
$$ |\ph_\cJ^{-1}(x)|=|\ker(\ph_\cJ)|=\frac{|G_\cJ|}{|\im(\ph_\cJ)|}
=\frac{\prod_{i\in\cJ}n_i}{n_\cJ} .$$
Thus the above expression for $|(\ph|_{\cO})^{-1}(x)|=m_d$ becomes
$$ (-1)^s\cdot m_d 
=\sum_{\cJ\sseq\cI,d\spmid n_\cJ}
\left(\frac{1}{n_\cJ}\cdot\prod_{i\in\cJ}(-n_i)\right) .$$

\Abs
This leads us to the second application of M\"obius inversion:
Letting
$$ k_c:=\sum_{\cJ\sseq\cI,n_\cJ=c}
\left(\frac{1}{c}\cdot\prod_{i\in\cJ}(-n_i)\right) $$
for all $c\spmid n$, we for all $d\spmid n$ get
\begin{equation*}\tag{$\ast$}
(-1)^s\cdot m_d=\sum_{c\in\{1\ld n\},d\spmid c\spmid n}k_c .
\end{equation*}
Hence considering this as a sum over the interval
$\{c\in\{1\ld n\};d\spmid c\spmid n\}$ 
with respect to the dual divisibility relation, M\"obius inversion yields
$$ (-1)^s\cdot k_d
=\sum_{c\in\{1\ld n\},d\spmid c\spmid n}\mu(\frac{c}{d})\cdot m_c ,$$
implying that $(-1)^s\cdot k_d$ can be determined from $\chi_\Phi$.

\abs 
Moreover, we have
$$ \sum_{c\in\{1\ld n\},c\spmid d}c\cdot k_c
=\sum_{c\in\{1\ld n\},c\spmid d}\left(\sum_{\cJ\sseq\cI,n_\cJ=c}
\left(\prod_{i\in\cJ}(-n_i)\right)\right) ,$$
yielding
$$ \sum_{c\in\{1\ld n\},c\spmid d}c\cdot k_c
=\sum_{\cJ\sseq\cI,n_\cJ\spmid d}\left(\prod_{i\in\cJ}(-n_i)\right) .$$
Letting $\cI_c:=\{i\in\cI,n_i=c\}$, for all $c\spmid n$,
the right hand side equals
$$ 
\sum_{\cJ\sseq\{i\in\cI,n_i\spmid d\}}\left(\prod_{i\in\cJ}(-n_i)\right)
=\prod_{i\in\cI,n_i\spmid d}(1-n_i)
=\prod_{c\in\{1\ld n\},c\spmid d}(1-c)^{|\cI_c|} .$$
Hence in conclusion we get, for all $d\spmid n$,
\begin{equation*}\tag{$\ast\ast$}
\sum_{c\in\{1\ld n\},c\spmid d}c\cdot k_c
=\prod_{c\in\{1\ld n\},c\spmid d}(1-c)^{|\cI_c|} .
\end{equation*}

\Abs\label{clarify}
We are now in a position to clarify the relationship
between the numbers $N:=\lcm(d\in\N;m_d>0)$ and $n_\cI=\lcm(n_i;i\in\cI)$.
As was remarked earlier, we have $m_d>0$ only if $d\spmid n_\cI$, 
thus we have $N\spmid n_\cI$. We consider the converse:

\abs
If $d\spnmid N$ then, by equation ($\ast$) we conclude that $k_d=0$. 
Hence equation ($\ast\ast$) applied to $n$ and $N$ yields
$$ \prod_{d\in\{1\ld n\},d\spmid n}(1-d)^{|\cI_d|}
=\sum_{d\spmid n}d\cdot k_d
=\sum_{d\spmid N}d\cdot k_d
=\prod_{d\in\{1\ld n\},d\spmid N}(1-d)^{|\cI_d|}, $$
showing that $|\cI_d|=0$ whenever $2\neq d\spnmid N$,
and $|\cI_2|$ is even whenever $N$ is odd.

\abs
Thus, if $N$ is even then 
$n_i\spmid N$ for all $i\in\cI$, 
hence $n_\cI\spmid N$, or equivalently $N=n_\cI$.
If $N$ is odd then we have two cases: If $n_\cI$ is odd as well,
then we have $n_i\spmid N$ for all $i\in\cI$, thus $N=n_\cI$ again;
if $n_\cI$ is even, then the multiset is of the form
$[n_1\ld n_{s'},2\ld 2]$, where $n_1\ld n_{s'}$ are odd and  
$s-s'=|\cI_2|>0$ is even, thus we have $2N=n_\cI$.

\abs
Hence, without knowing the underlying multiset $[n_1\ld n_s]$,
we are able to specify a small multiple $n$ of $n_{\cI}$,
using $\chi_\Phi$ alone, as
$$ n=\left\{\begin{array}{rl}
N,&\text{ if $N$ is even},\\
2N,&\text{ if $N$ is odd}.\\
\end{array} \right. $$
Moreover, assuming that $[n_1\ld n_s]$ contains the element $2$
at most once, then the above also shows that the order of 
periodicity of $\Phi$ is given as $N=n_\cI$,
showing the first half of Corollary \ref{cor1}.
 
\Abs
We are now prepared to describe how to recover the multiset $[n_1\ld n_s]$ 
from $\chi_\Phi$:
Firstly, for $d=1$ we have $|\cI_1|=0$, 
thus by equation ($\ast\ast$) we get $k_1=0^0=1$,
hence being able to compute $(-1)^s\cdot k_1=(-1)^s$ 
implies that the parity of the number $s$
of tensor factors can be determined from $\chi_\Phi$.
Thus we are able to compute the numbers $k_d$, for all $d\spmid n$.

\abs
Next, for $d=2$ equation ($\ast\ast$) yields $1+2k_2=(-1)^{|\cI_2|}$, thus 
$k_2=\frac{(-1)^{|\cI_2|}-1}{2}\in\{0,-1\}$,
where $k_2=0$ if and only if $|\cI_2|$ is even.
Hence being able to compute $k_2$
implies that the parity of the number $|\cI_2|$ of tensor factors 
equal to $2$ can be determined from $\chi_\Phi$ as well.
Recall that we assume that $n$ is even, so that the above
analysis indeed holds for $d=2$, and that 
since $\chi_{[n_1\ld n_s]}=\chi_{[n_1\ld n_s,2,2]}$
we cannot possibly expect more than finding the parity of $|\cI_2|$.

\abs
Finally, for divisors $d\spmid n$ such that $d>2$ we proceed by induction,
using equation ($\ast\ast$) again, to obtain successively $(1-d)^{|\cI_d|}$, 
which since $d-1>1$ immediately yields $|\cI_d|$.

\abs
In conclusion, for all $d\spmid n$, 
we are able to recover algorithmically the cardinality $|\cI_d|$
of the number of tensor factors equal to $d$ if $d>2$,
and the parity $|\cI_2|$ of the number of tensor factors 
equal to $2$.
This in particular proves Theorem \ref{thm2}.
\QED

\Ex\label{ex2}
For example, for the multisets $[3,4,5,6,7]$ and $[2,3,4,5,6,7]$ 
we get $\deg(\chi_\Phi)=6\cdot 5\cdot 4\cdot 3\cdot 2=720$,
and the following non-zero multiplicities 
$m_d$ and $m'_d$ of zeroes of $\chi_\Phi$, respectively,
entailing $n=N=n_\cI=420$:
$$ \begin{array}{|r||r|r|r|r|r|r|}
\hline 
d&35&70&105&140&210&420\\
\hline 
\hline 
m_d&2&2&2&4&1&3\\
\hline 
m'_d&2&2&1&4&2&3\\
\hline 
\end{array} $$


\abs
This yields the associated numbers $k_d$ and $k'_d$, respectively,
as follows:
$$ \begin{array}{|r||r|r|r|r|r|r|r|r|r|r|r|r|r|}
\hline 
d&1&2&3&4&5&6&7&10&12&14&15&20&21\\
\hline 
\hline 
(-1)^s\cdot k_d&-1&0&1&1&1&-2&1&0&3&0&-1&-1&-1\\
\hline 
(-1)^s\cdot k'_d&1&-1&-1&1&-1&-1&-1&1&3&1&1&-1&1\\
\hline 
\end{array} $$
$$ \begin{array}{|r||r|r|r|r|r|r|r|r|r|r|r|}
\hline 
d&28&30&35&42&60&70&84&105&140&210&420\\
\hline 
\hline 
(-1)^s\cdot k_d&-1&2&-1&2&-3&0&-3&1&1&-2&3\\
\hline 
(-1)^s\cdot k'_d&-1&1&1&1&-3&-1&-3&-1&1&-1&3\\
\hline 
\end{array} $$

\abs
These now yield, by considering $d=1$, that $s$ is odd respectively even, 
and by considering $d=2$, that $|\cI_2|$ is even respectively odd,
and then by induction on $d$ reveal successively
$|\cI_3|=|\cI_4|=|\cI_5|=|\cI_6|=|\cI_7|=1$, and $|\cI_d|=0$
for $d\spmid n$ such that $d\geq 8$.

\section{Spectra of Coxeter transformations}\label{coxspec}

\abs
We finally comment on two aspects concerning explicit eigenvalues
of Coxeter transformations are concerned, thereby proving the second
half of Corollary \ref{cor1}, and Corollary \ref{cor2}.
From now on we assume that the multiset $[n_1\ld n_s]$
contains the element $2$ at most once.

\AbsT{The root of unity $\zt_{n_\cI}$ as an eigenvalue.}\label{primedec}
In order to prove the second half of Corollary \ref{cor1}, 
we proceed to show that $m_{n_\cI}>0$, in particular implying 
$$ \max(d\in\N;m_d>0)=\lcm(d\in\N;m_d>0)=N .$$
Now, asking whether $m_{n_\cI}>0$ 
is equivalent to asking whether the image $\im(\ph|_\cO)\sseq\Z/n$ 
contains an element of additive order $n_\cI$. Since the map
$$ \ph=[\ia_n^{n_1}\ld\ia_n^{n_s}]\cn\Z/n_1\oplus\cdots\oplus\Z/n_s\ra\Z/n $$
factors through $\ia_n^{n_\cI}\cn\Z/n_\cI\ra\Z/n$ anyway, 
we may assume that $n=n_\cI$, and ask whether $1\in\im(\ph|_\cO)\sseq\Z/n$.

\abs
To this end, let $p_1\ld p_t\in\N$ be the rational prime divisors
of $n$, for some $t\in\N$, and let $q_j:=p_j^{\eps_j}$, where $\eps_j\in\N$,
such that $n=\prod_{j=1}^t q_j$.
For $i\in\{1\ld s\}$ and $j\in\{1\ld t\}$ let 
$q_{ij}:=p_j^{e_{ij}}$, where $e_{ij}\in\N_0$, such that
$n_i=\prod_{j=1}^t q_{ij}$.
Thus we have the primary decompositions 
$$ \Z/n_i\cong\Z/q_{i1}\oplus\cdots\oplus\Z/q_{it}
\quad\text{ and }\quad
\Z/n\cong\Z/q_1\oplus\cdots\oplus\Z/q_t ,$$
given by the associated natural embeddings.
Note that for those we have the following commutative diagram:
$$\begin{CD}
\Z/n_i @>{\ia_n^{n_i}}>> \Z/n\\
@A{\ia_{n_i}^{q_{ij}}}AA @AA{\ia_n^{q_j}}A\\
\Z/q_{ij} @>>{\ia_{q_j}^{q_{ij}}}> \Z/q_j \\
\end{CD}$$
Hence for $j\in\{1\ld t\}$ letting
$$ \ph_j\cn\bigoplus_{i\in\{1\ld s\};e_{ij}>0}\Z/q_{ij}\ra\Z/q_j ,$$
we conclude that $1\in\im(\ph)\sseq\Z/n$, that is
$\im(\ph)$ contains an element of additive order $n$,
if and only if $\im(\ph_j)$ contains an element of additive order $q_j$,
that is $1\in\im(\ph_j)\sseq\Z/q_j$, for all $j\in\{1\ld t\}$.
We consider restrictions to $\cO$:

\abs
Assume that $1\in\im(\ph_j|_\cO)\sseq\Z/q_j$, for all $j\in\{1\ld t\}$.
Taking a preimage in $\ph_j^{-1}(1)\cap\cO$,
and augmenting by $0\in\Z/q_{ij}$ for $i\in \{1\ld s\}$ such
that $e_{ij}=0$, we obtain an element
$x_j=[x_{1j}\ld x_{sj}]\in\Z/q_{1j}\oplus\cdots\oplus\Z/q_{sj}$.
Letting 
$$ y_{ij}:=\ia_{n_i}^{q_{ij}}(x_{ij})\in\im(\ia_{n_i}^{q_{ij}})\leq\Z/n_i ,$$
for $i\in \{1\ld s\}$, we have $y_{ij}\neq 0$ if and only if $e_{ij}>0$,
where in this case $y_{ij}$ has additive order $q_{ij}$, and 
\begin{align*}
\ph([y_{1j}\ld y_{sj}])
&=\sum_{i=1}^s\ia_n^{n_i}(\ia_{n_i}^{q_{ij}}(x_{ij}))
 =\sum_{i=1}^s\ia_n^{q_j}(\ia_{q_j}^{q_{ij}}(x_{ij})) \\
&=\ia_n^{q_j}\left(\sum_{i=1}^s\ia_{q_j}^{q_{ij}}(x_{ij})\right) 
 =\ia_n^{q_j}(\ph_j(x_j))
 =\ia_n^{q_j}(1) 
 =\frac{n}{q_j} 
\end{align*}
shows that $\ph([y_{1j}\ld y_{sj}])\in\im(\ia_n^{q_j})\leq\Z/n$ 
has additive order $q_j$.
Then for 
$$ y:=\sum_{j=1}^t\, [y_{1j}\ld y_{sj}]
=[\sum_{j=1}^t y_{1j}\ld\sum_{j=1}^t y_{sj}]
\in\Z/n_1\oplus\cdots\oplus\Z/n_s
$$
we have $\sum_{j=1}^t y_{ij}\neq 0\in\Z/n_i$,
that is $y\in\cO$, and $\ph(y)=\sum_{j=1}^t \frac{n}{q_j}\in\Z/n$
has additive order $\lcm(q_1\ld q_s)=n$,
hence we indeed have $1\in\im(\ph|_\cO)\sseq\Z/n$.

\Abs\label{oneprime}
Thus to show that $1\in\im(\ph|_\cO)\sseq\Z/n$ 
we are led to consider the case $t=1$:
To simplify notation, 
let $p\in\N$ be a rational prime, and let $n=p^e$ for some $e\in\N$.
Hence we have $n_i=p^{e_i}$, where $e_i\in\N$, and where
we may assume that $e_i<e$ for all $i\in\{1\ld s'\}$,
and $e_i=e$ for all $i\in\{s'+1\ld s\}$, for some $s'\in\{0\ld s-1\}$.
Now we have $1\in\im(\ph|_\cO)\sseq\Z/p^e$ if and only if $m_{p^e}\neq 0$,
which, by equation $(\ast)$ applied to $d=p^e$, 
is equivalent to $k_{p^e}\neq 0$.

\abs 
To determine $k_{p^e}$, for $\cJ\sseq\cI$
we have $n_{\cJ}=p^e$ if and only if $\cJ\cap\{s'+1\ld s\}\neq\emp$,
that is $\cJ$ is of the form $\cJ=\cJ'\dcup\cJ''$ where
$\cJ'\sseq\{1\ld s'\}$ and $\emp\neq\cJ''\sseq\{s'+1\ld s\}$.
Thus the principle of inclusion-exclusion yields
\begin{align*}
p^e\cdot k_{p^e}&=
\sum_{\cJ\sseq\cI,\cJ\cap\{s'+1\ld s\}\neq\emp}
\left(\prod_{i\in\cJ}(-p^{e_i})\right) \\
&=
\sum_{\cJ'\sseq\{1\ld s'\}}\,\,\,\sum_{\emp\neq\cJ''\sseq\{s'+1\ld s\}}
\left(\prod_{i\in\cJ'\dcup\cJ''}(-p^{e_i})\right) \\
&=
\sum_{\cJ'\sseq\{1\ld s'\}}\left(\prod_{i\in\cJ'}(-p^{e_i})\right)
\,\,\,\cdot\,\,\,
\sum_{\emp\neq\cJ''\sseq\{s'+1\ld s\}}\left(\prod_{i\in\cJ''}(-p^e)\right)\\
&=
\left(\prod_{i\in\{1\ld s'\}}(1-p^{e_i})\right)\cdot
\left((1-p^e)^{s-s'}-1\right) .
\end{align*}
Hence we conclude that $k_{p^e}\neq 0$, except if 
$p^e=2$ and $s-s'=s$ is even.

\Abs
Thus, returning to the general case again, by the above reduction 
we are done in all cases, except if $n$ is even and $\eps_1=1$,
where we let $p_1:=2$. In this case, we may assume that there is 
$s'\in\{0\ld s-1\}$ such that $n_1\ld n_{s'}$ 
are odd, and $n_{s'+1}\ld n_s$ are even. By the above reduction again
we are done if $s-s'>0$ is odd,
hence we may additionally assume that $s-s'\geq 2$ is even. 

\abs
To simplify notation, let $n'_i:=n_i$ for $i\in\{1\ld s'\}$,
and $n'_i:=\frac{n_i}{2}\in\N$ for $i\in\{s'+1\ld s\}$, and let
$n':=\frac{n}{2}=\lcm(n'_1\ld n'_s)$.
Hence the $n'_i$ are odd, where $n'_i=1$ if and only if $n_i=2$.
Since the multiset $[n_1\ld n_s]$ contains the element $2$ at most once, 
we may assume that $n'_i\geq 3$ for all $i\in\{1\ld s-1\}$, 
and we have to distinguish the cases $n'_s=1$ and $n'_s\geq 3$:

\abs
If $n'_s\geq 3$, we consider the map
$\ph'\cn\Z/n'_1\oplus\cdots\oplus\Z/n'_s\ra\Z/n'$.
Then by the above analysis we have $1\in\im(\ph'|_\cO)\sseq\Z/n'$,
hence there is 
$$ [y_1\ld y_s]\in\ia_{n_1}^{n'_1}((\Z/n'_1)\smin\{0\})
   \oplus\cdots\oplus\ia_{n_s}^{n'_s}((\Z/n'_s)\smin\{0\})\sseq\cO $$
such that $\ph([y_1\ld y_s])\in\Z/n$ has additive order $n'$.
Since $y_s\in\ia_{n_s}^{n'_s}(\Z/n'_s)\sseq\Z/n_s$ has odd order,
and $n'_s=\frac{n_s}{2}\in\Z/n_s$ has additive order $2$,
we conclude that $y_s+n'_s\in(\Z/n_s)\smin\{0\}$ has even order.
Hence we have $[y_1\ld y_{s-1},y_s+n'_s]\in\cO$, and 
$\ph([y_1\ld y_{s-1},y_s+n'_s])\in\Z/n$
has additive order $n=2n'$.

\abs
If $n'_s=1$, we consider the map
$\ph'\cn\Z/n'_1\oplus\cdots\oplus\Z/n'_{s-1}\ra\Z/n'$ instead.
Then similarly there is 
$$ [y_1\ld y_{s-1}]\in\ia_{n_1}^{n'_1}((\Z/n'_1)\smin\{0\})
   \oplus\cdots\oplus 
   \ia_{n_{s-1}}^{n'_{s-1}}((\Z/n'_{s-1})\smin\{0\}) $$
such that $\ph([y_1\ld y_{s-1},0])\in\Z/n$ has additive order $n'$.
Thus, since $1\in(\Z/n_s)\smin\{0\}=\{1\}$ has has additive order $2$,
we again have $[y_1\ld y_{s-1},1]\in\cO$, and
$\ph([y_1\ld y_{s-1},1])\in\Z/n$ has additive order $n=2n'$.
\QED

\AbsT{The root of unity $1$ as an eigenvalue.}\label{primeset}
Quite to the opposite we now deal with the question whether $m_1>0$,
or equivalently whether $0\in\im(\ph|_\cO)\sseq\Z/n$.
We proceed similar to \ref{primedec}:
If $0\in\im(\ph_j|_\cO)\sseq\Z/q_j$, for some $j\in\{1\ld t\}$,
there again are $y_{ij}\in\im(\ia_{n_i}^{q_{ij}})\leq\Z/n_i$, 
for $i\in \{1\ld s\}$, 
such that $y_{ij}\neq 0$ if and only if $e_{ij}>0$,
where in this case $y_{ij}$ has additive order $q_{ij}$,
and $\ph([y_{1j}\ld y_{sj}])=0\in\Z/n$.

\abs
Thus to infer $0\in\im(\ph|_\cO)\sseq\Z/n$ 
it suffices to assume that $0\in\im(\ph_j|_\cO)\sseq\Z/q_j$ 
where $j$ runs through a subset $\cK\sseq\{1\ld t\}$ such that 
$\sum_{j\in\cK}[y_{1j}\ld y_{sj}]\in\cO$.
The latter condition, saying
$\sum_{j\in\cK}y_{ij}\neq 0\in\Z/n_i$ for all $i\in\{1\ld s\}$,
is equivalent to
$\{j\in\{1\ld t\};e_{ij}>0\}\cap\cK\neq\emp$, for all $i\in\{1\ld s\}$,
in other words $\{p_j;j\in\cK\}$ 
contains a prime divisor of any of the numbers $n_1\ld n_s$.

\Abs\label{oneprimeagain}
Thus to show that $0\in\im(\ph|_\cO)\sseq\Z/n$
we again first consider the case $t=1$:
We keep the notation of \ref{oneprime},
where we also assume that $e_1\leq\cdots\leq e_s\in\N$.
Then by equation $(\ast)$, applied to $d=1$, we have
$$ (-1)^s\cdot m_1
=\sum_{\cJ\sseq\cI}
\left(\frac{1}{n_\cJ}\cdot\prod_{i\in\cJ}(-p^{e_i})\right) 
=1+\sum_{\emp\neq\cJ\sseq\cI}
\left(\frac{1}{p^{e_{\max(\cJ)}}}\cdot\prod_{i\in\cJ}(-p^{e_i})\right) .$$
Reordering the sum with respect to the maximum of the
subsets $\emp\neq\cJ\sseq\cI$ considered, the right hand side becomes
\begin{align*}
&1+\sum_{\emp\neq\cJ\sseq\cI}
\left(\frac{1}{p^{e_{\max(\cJ)}}}\cdot\prod_{i\in\cJ}(-p^{e_i})\right) 
=1-\sum_{\emp\neq\cJ\sseq\cI}
\left(\prod_{i\in\cJ\smin\max(\cJ)}(-p^{e_i})\right) \\
=\,\,&1-\sum_{j=1}^s\,\,\,\sum_{\cJ\sseq\{1\ld j-1\}}
\left(\prod_{i\in\cJ}(-p^{e_i})\right) 
=1-\sum_{j=1}^s\left(\prod_{i=1}^{j-1}(1-p^{e_i})\right) .
\end{align*}
Hence we conclude that
$$ m_1=\sum_{j=2}^s \left((-1)^{s-j}\cdot
\prod_{i=1}^{j-1}(p^{e_i}-1)\right) .$$

\abs
Thus we have $m_1=0$ whenever $s=1$. Hence let $s\geq 2$,
then $m_1$ is an alternating sum with 
summands having non-decreasing absolute value,
even increasing absolute value for those $j$ such that $p^{e_{j-1}}\geq 3$.
Hence we conclude that $m_1>0$, except if $p=2$ and 
$1=e_1=\cdots=e_{s-1}\leq e_s=e$, where $s$ is odd.

\Abs
Now, returning to the general case again,
to state the criterion to decide whether $m_1>0$
we need a few preparations, taking the above analysis into account:
Let $\Gm$ be the {\bf gcd graph} of the multiset $[n_1\ld n_s]$,
whose vertices are labelled by $n_1\ld n_s$, and where vertices
$n_i$ and $n_{i'}$, for $i\neq i'\in\{1\ld s\}$, 
are adjacent if and only if $\gcd(n_i,n_{i'})>1$.

\abs
The graph $\Gm$ is the union of the various {\bf prime graphs} 
$\Gm_j$ of $[n_1\ld n_s]$, for $j\in\{1\ld t\}$, 
which also have vertices labelled by $n_1\ld n_s$,
where vertices $n_i$ and $n_{i'}$, for $i\neq i'\in\{1\ld s\}$, 
are adjacent in $\Gm_j$ if and only if $p_j\spmid\gcd(n_i,n_{i'})$. 
Hence each of the graphs $\Gm_j$, next to isolated vertices, 
has at most one non-trivial connected component $\Gm_j^0$. If it exists,
$\Gm_j^0$ is a complete graph with at least two vertices, being labelled
by the $n_i$ having $j$-th exponent $e_{ij}>0$; if it does not exist
we for completeness let $\Gm_j^0$ be the empty graph.

\abs
In particular, $n_i$ is an isolated vertex of $\Gm$ if and only if 
$n_i$ does not belong to any of the connected components $\Gm_j^0$.
Moreover, if $n_i$ is an isolated vertex of $\Gm_j$ then we are
in the case $s=1$ in the analysis in \ref{oneprimeagain}, 
thus we have $0\not\in\im(\ph_j|_\cO)$,
that is all elements of $\im(\ph_j|_\cO)\sseq\Z/q_j$
have additive order divisible by $p_j$. Thus, in this case, if
$z=[z_1\ld z_s]\in\ph^{-1}(0)\sseq\Z/n_1\oplus\cdots\oplus\Z/n_s$,
then writing 
$z_i=[z_{i1}\ld z_{it}]\in\Z/q_{i1}\oplus\cdots\oplus\Z/q_{it}\cong\Z/n_i$
we conclude that $z_{ij}=0\in\Z/q_{ij}$.

\Thm\label{m1thm}
Let $[n_1\ld n_s]$ be a multiset containing the element $2$ at most once.
Then $1$ is \emph{not} an eigenvalue of $\Phi_{[n_1\ld n_s]}$ 
if and only if 

\abs\bfi 
$\Gm$ has an isolated vertex, or 

\abs\bfii 
$n_\cI$ is even, and letting $p_1:=2$ the graph
$\Gm_1^0$ has an odd number of vertices, none of which belongs to any
of the $\Gm_j^0$, where $j\in\{2\ld t\}$, and the associated
multiset of exponents is $[1\ld 1,\eps]$ for some $\eps\geq 1$.

\Pf
We first show that $m_1=0$ if either of conditions (i) or (ii) holds:
Firstly, if $\Gm$ has an isolated vertex $n_i$, then it is 
isolated in $\Gm_j$ as well, for all $j\in\{1\ld t\}$.
Hence for any $z\in\ph^{-1}(0)$, keeping the above notation,
we have $z_{ij}=0\in\Z/q_{ij}$ for all $j\in\{1\ld t\}$,
that is $z_i=0\in\Z/n_i$, hence $z\not\in\cO$.


\abs
Secondly, we may assume that $n_1\ld n_{s'}$ are odd, 
and $n_{s'+1}\ld n_s$ are even, for some $s'\in\{0\ld s-1\}$.
Then the vertices of $\Gm_1^0$ are labelled by $[n_{s'+1}\ld n_s]$, 
hence $s-s'>0$ is odd, and since we may assume that $\Gm$ does not
have isolated vertices, we may additionally assume that $s-s'\geq 3$.
Letting $i\in\{s'+1\ld s\}$, since $n_i$ is an 
isolated vertex of $\Gm_j$ for all $j\in\{2\ld t\}$, 
we may assume that the vertices of 
$\Gm_1^0$ are labelled by the multiset $[2\ld 2,2^\eps]$.
Hence for any preimage $z=[z_1\ld z_s]\in\ph^{-1}(0)\cap\cO$ we infer
$[z_{s'+1}\ld z_s]\in\ph_1^{-1}(0)\cap\cO$,
which is the exceptional case in \ref{oneprimeagain}, a contradiction.

\abs
Now assume that neither of conditions (i) and (ii) hold.
Then we aim to show that the set 
$\cK:=\{j\in\{1\ld t\};\Gm_j^0\text{ non-empty}\}$
fulfils the conditions described in \ref{primeset}.
Note that, since $\Gm$ does not have isolated vertices, 
the set $\{p_j;j\in\cK\}$ contains a prime divisor of
any of the numbers $n_1\ld n_s$. We are going to apply \ref{oneprimeagain}
repeatedly: First of all, we are done if $n_\cI$ is odd.
Hence we may assume that $n_\cI$ is even, and let $p_1:=2$, 
thus $q_1=2^\eps$ for some $\eps\in\N$. Then, if $\Gm_1^0$ is empty, 
or is non-empty and has an even number of vertices, we are done as well.

\abs
Hence we may assume that $\Gm_1^0$ has an odd number of vertices,
being labelled by $[n_{s'+1}\ld n_s]$, for some $s'\in\{0\ld s-1\}$ 
such that $s-s'\geq 3$ is odd.
Assume next that $n_s$, say, belongs to $\Gm_1^0$ as well as to $\Gm_j^0$
for some $j\in\{2\ld t\}$. Then we replace the map $\ph_1$ by
$$ \ph'_1\cn\Z/q_{s'+1,1}\oplus\cdots\oplus\Z/q_{s-1,1}\ra\Z/q_1=\Z/2^\eps ,$$
where the left hand side has an even number of summands.
Hence we have $0\in\im(\ph'_1|_\cO)$, and the conditions in 
\ref{primeset} hold for $\cK$ and the maps $\ph'_1,\ph_2\ld\ph_t$.

\abs
Finally, assume that none of the vertices of $\Gm_1^0$ 
belong to any of the $\Gm_j^0$, where $j\in\{2\ld t\}$.
Thus by assumption the $2$-parts $[q_{s'+1,1}\ld q_{s1}]$ 
are not of the form $[2\ld 2,2^\eps]$, implying that
$0\in\im(\ph_1|_\cO)\sseq\Z/2^\eps$, and we are done.
\QED

\Ex\label{ex3}
Since the multisets given \ref{ex2} contain isolated vertices,
entailing $m_1=0$, we here give a few examples without isolated vertices:
For the multisets $[2,4,6]$ and $[2,3,4,6]$ and $[2,4,6,6]$
the associated characteristic polynomials have degree
$15$, $30$ and $75$, respectively, and we get the following associated
multiplicities $m_d$, $m'_d$ and  $m''_d$, respectively:
$$ \begin{array}{|r||r|r|r|r|r|r|}
\hline
d&1&2&3&4&6&12\\
\hline
\hline
m_d&0&1&1&1&1&2\\
\hline
m'_d&2&2&1&4&2&3\\
\hline
m''_d&5&4&4&9&4&8\\
\hline
\end{array} $$

\Abs
We finally prove Corollary \ref{cor2}, where it indeed suffices to
assume that the multiset $[n_1\ld n_s]$
contains the element $2$ at most once:
The process of taking reciprocal polynomials is
an involutory automorphism of the multiplicative semigroup 
$\Z[X]\smin\{0\}$. Moreover, letting $\Phi_d\in\Z[X]$ be the
$d$-th cyclotomic polynomial, where $d\in\N$, then we have
$\Phi_1^\ast=(X-1)^\ast=-(X-1)=-\Phi_1$, while $\Phi_d$
is self-reciprocal for $d\geq 2$,
see for example \cite[Sect.1.3]{LenzingPenaII}.
Since $\chi_{[n_1\ld n_s]}\in\Z[X]$ is a product of
cyclotomic polynomials, we thus conclude that 
$\chi_{[n_1\ld n_s]}^\ast=\pm\chi_{[n_1\ld n_s]}$.

\abs
Hence for the constant term of $\chi_{[n_1\ld n_s]}$ we have
$\chi_{[n_1\ld n_s]}(0)=\pm 1$, where $\chi_{[n_1\ld n_s]}$ is 
self-reciprocal if and only if $\chi_{[n_1\ld n_s]}(0)=1$. Thus from 
$$ \chi_{[n_1\ld n_s]}(0)=(-1)^l\cdot\det(\Phi_{[n_1\ld n_s]})
   =(-1)^{l(s+1)} ,$$
where $l:=\prod_{i=1}^s(n_i-1)$,
we infer that $\chi_{[n_1\ld n_s]}$ is self-reciprocal if and only if
$s$ is odd or at least one of the $n_i$ is odd.

\abs
(Alternatively, this also follows from recalling that 
$\chi_{[n_1\ld n_s]}$ or $\chi_{[n_1\ld n_s,2]}$ is a Coxeter polynomial, 
depending on whether $s$ is odd or even, respectively,
and using the fact that Coxeter polynomials always are self-reciprocal.) 
\QED 

\abs\abs
Finally, note that by the above analysis 
$\chi_{[n_1\ld n_s]}$ is self-reciprocal
if and only if the multiplicity $m_1\in\N_0$ of the
eigenvalue $1$ is even; see also the examples in \ref{ex3}.
Hence the conditions characterising the case $m_1=0$ 
in Theorem \ref{m1thm} constitute a special case of
the condition in Corollary \ref{cor2}: Indeed, 
the conditions in Theorem \ref{m1thm} imply that $\Gm$ is not connected,
which in turn entails that $\gcd(n_1\ld n_s)=1$, truly a 
special case of $\gcd(s,n_1\ld n_s)$ being odd.


\abs
{\sc L.H.: 
Mathematisches Institut, Universit\"at M\"unster \\
Einsteinstra{\ss}e 62, D-48149 M\"unster, Germany} \\
{\sf lutz.hille@uni-muenster.de }

\abs
{\sc J.M.: 
Lehrstuhl D f\"ur Mathematik, RWTH Aachen \\
Templergraben 64, D-52062 Aachen, Germany} \\
{\sf juergen.jueller@math.rwth-aachen.de}

\end{document}